\documentclass[a4paper]{article} 

\newtheorem{algthm}{Algorithm}[section]
\newcounter{alg}

\begin{document}











\title{Alternative Local Discriminant Bases Using Empirical
  Expectation and Variance Estimation}

\author{Eirik Fossgaard,\thanks{Eirik Fossgaard is currently a
    lecturer at Institute of Mathematics and Statistics, University of
    Troms{\o}, 9037 Troms{\o}, Norway (E-mail: eirikf@math.uit.no). 
    This work was supported by Department of Mathematics, Royal Institute of Technology, 100 44
    Stockholm, Sweden. The author thanks Professor Jan Olov
    Str\"{o}mberg, Department of Mathematics, Royal Institute of
    Technology, 100 44 Stockholm, Sweden, for helpful suggestions.}}

\maketitle

\begin{abstract}
We propose alternative discriminant measures for selecting the best basis among a
large collection of orthonormal bases for classification purposes. A
generalization of the Local Discriminant Basis Algorithm of Saito and
Coifman is constructed. The success of these new methods is evaluated
and compared to earlier methods in experiments.
\end{abstract}



\section{Introduction}

This paper is the result of my trying to improve the method
applied in Fossgaard (1997) to discriminate between two distinct
classes/types of signals by using expansions of the data in 
wavelet packet/local trigonometric bases. This method was first
invented and described by N.Saito and R.Coifman. For a thorough
exposition on this theme, I refer to Saito (1994) and Saito, Coifman
(1996), a brief summary of the main ideas is given below.

Each signal belonging to a training dataset is decomposed in a 
time/space -frequency {\em dictionary}, that is a decomposition into a
large collection of orthonormal bases arranged in a binary tree structure, containing 
either {\em wavelet-packet} basis functions, or 
{\em local trigonometric} basis functions. A measure of {\em energy-density}
is then computed for each coordinate in the dictionary for each class
of signals, originally in Saito (1994) this is taken
to be the square of the coordinate summed over all the training signals
belonging to a class of signals, and then 
normalized by the total energy projected onto this coordinate. 
Then a basis called the ``Local Discriminant Basis'', LDB for
short, is chosen from the dictionary by maximizing a certain {\em
  discrimination measure}, defined by some {\em additive} cost-functional, 
over the dictionaries of energy-densities. 
The coordinates where the discrimination measure takes on its largest
values are called the most important {\em features} of the signals.
These coordinates are selected from the LDB and used as input for some classifier.

This method is very powerful in many cases, but it also has its
weaknesses, a serious one is that the LDB is not able to distinguish 
two signals both consisting exclusively of one and the same basis
element, only with opposite sign. One way of dealing with this problem
is described in Saito, Coifman (1996), where one estimates the {\em
  probability-density functions}, pdf's, of the projections onto the different basis elements in
the dictionary, and selects the basis which maximizes some well-chosen
functional on these pdf's. 

In this paper, I will try to improve on the LDB-method described above, 
by constructing new dicrimination measures that yield more relevant 
features. I will also try to improve the 
performance of the algorithm by using {\em several}
LDB's in sequence, and by using a classifier specially designed to fully utilize   
the increased degree of freedom multiple LDB's (MLDB's) give us in
selecting features that are most important to our problem.

\section{The original LDB method}
 
The problem as expressed in Saito (1994) 
is optimizing a linear map: $d:\mathcal{X}\rightarrow \mathcal{Y}$, 
where $\cup_{y\in{\mathcal Y}}{\mathcal X}^{(y)} = {\mathcal X}\subset
{\mathbf R}^{n}$  is the input signal space,
${\mathcal Y}=\{1,2,...,N\}$ is the output class space, a set of class
labels, and ${\mathcal X}^{(y)}$ is the subspace of class $y$ signals.
To optimize the map $d$, one considers maps
of the form  
\begin{equation}  
d=c\circ{\mathcal F}_{K}\circ{\mathbf \Psi}_{n\times n}, 
\label{map_problem}
\end{equation}
where the {\em feature extractor} ${\mathbf \Psi}_{n\times n}\in O(n)$ is an
orthogonal $n\times n$ matrix which extracts the $n$ most relevant
coordinates from from a binary-tree dictionary 
of wavelet packet bases or local trigonometric bases, ${\mathcal F}_{K}$
is a {\em feature selector} which selects the $K<n$ most 
important coordinates from the $n$ most relevant coordinates, and $c$ is a
classifier. The problem then is to choose $c,{\mathcal F}_{K}$ and
${\mathbf \Psi}_{n\times n}$ such that the rate of misclassification of the map 
$d$ is minimized on the set ${\mathcal X}$. 
In Saito (1994), ${\mathbf \Psi}_{n\times n}$ is taken to be
\begin{equation}
  {\mathbf \Psi}_{n \times n} = \mbox{arg}\max_{B_{k}\in
    {\mathcal D}_{i}\in{\mathcal L}}\lambda(B_{k}),
  \label{best_basis}
\end{equation} 
where ${\mathcal L}=\cup_{i}{\mathcal D}_{i}$  is the {\em library} of  
all dictionaries at our disposal   
corresponding to the different wavelet or local trigonometric
basis functions under consideration, the $B_{k}$ are all bases in
${\mathcal D}_{i}$, and $\lambda$ is a measure of performance of
the basis $B_{k}$ in the classification problem, such a measure is
called a discrimination measure. The search for this
${\mathbf \Psi}_{n\times n}$ is fast by the {\em best-basis-algorithm} of Wickerhauser
and Coifman if the measure $\lambda$ satisfies an additivity property, 
(Saito 1994). In Saito (1994) the discrimination measure $\lambda$ is defined as
\begin{equation}
  \lambda(B_{k})=\sum_{{\mathbf w}_{m}\in B_{k}}
  \gamma(\Gamma^{(1)}({\mathbf w}_{m}),...,\Gamma^{(N)}({\mathbf w}_{m})),
  \label{performance_measure}
\end{equation}
where the {\em time-frequency energy-map} $\Gamma^{(y)}$ is defined by
\begin{eqnarray}
  && \Gamma^{(y)}({\mathbf w}_{m}) = 
  \frac{\sum_{j=1}^{J_{y}}({\mathbf w}_{m}\cdot{\mathbf x}^{(y)}_{j})^{2}}
  {\sum_{j=1}^{J_{y}}\|{\mathbf x}^{(y)}_{j}\|^{2}},\nonumber \\
  && {\mathbf x}^{(y)}_{j}\in{\mathcal X}^{(y)},\mbox{
    }1\leq y\leq N, \ J_{y}=|{\mathcal X}^{(y)}|,  
  \label{Gamma}
\end{eqnarray}
and $\gamma$ can be some form of $l^{p}$- distance, Hellinger-distance 
or relative entropy. \\

The signals ${\mathbf x}^{(y)}\in {\mathcal X}^{(y)}, \mbox{  }
1\leq y \leq N$ are fed into each dictionary as given by
(\ref{Gamma}), the best basis picked
out by the best basis algorithm, and then the best $K$ coordinates are
selected from this basis, ordinarily by selecting the coordinates
where $\lambda$ takes on its $K$ greatest values. The corresponding $K$ best basis elements
are then used to construct a classifier by doing a ``Linear Discriminant
Analysis'' (LDA) or a ``Classification and Regression Trees''
(CART)-analysis, or some other statistical classification technique, on the coordinates 
of the signals in these $K$ best basis elements. 

\section{A generalized LDB method}
\subsection{New Discrimination Measures}
Using the notation from the previous section, for each basis vector
${\mathbf w}_{m}$ in some basis $B_{k}$, let $Z_{y,m}$ be random
variable on the space ${\mathcal X}^{(y)}$ of input signals of class
$y$ defined by 
\begin{equation}
  Z_{y,m}: {\mathbf x}\in{\mathcal X}^{(y)}\rightarrow [-1,1],\
  Z_{y,m}({\mathbf x}) = {\mathbf w}_{m}\cdot{\mathbf x}.
  \label{random_Z}
\end{equation}

\noindent In Saito, Coifman (1996) one estimates the empirical pdf $p$
of $Z_{y,m}$. These estimates are then used to find the most
discriminating basis. But getting good estimates of the pdf's is hard
and computationally demanding. We will take a different approach and
work on  the {\em a priori} assumption that $p$ is the uniform distribution.
For each fixed ${\mathbf w}_{m}\in B_{k}$, we can then compute
the {\em empirical expectation} $E[Z_{y,m}]$
of the basis coordinate ${\mathbf w}_{m}\cdot{\mathbf x}$ for class $y$ signals as 
\begin{eqnarray}
  E[Z_{y,m}] &=& \sum_{{\mathcal X}^{(y)}}p(Z_{y,m}|Y=y)Z_{y,m}\nonumber\\
  &=&\sum_{{\mathbf x}\in{\mathcal X}^{(y)}}\frac{1}{|{\mathcal X}^{(y)}|}({\mathbf
    w}_{m}\cdot{\mathbf x}).
\end{eqnarray}
If $\|{\mathbf x}\|_{2} = 1, \forall {\mathbf x}\in {\mathcal X}$, then in this
probabilistic setting, (\ref{Gamma}) is equivalent to
$\Gamma^{(y)}({\mathbf w}_{m}) = E[Z^{2}_{y,m}]$. We will first
consider two-class problems: ${\mathcal Y}=\{1,2\}$, and deal with
$n$-class problems later. Choosing $\gamma =
\ell^{2}-distance\mbox{ }squared$, (\ref{performance_measure}) becomes
\begin{equation}
\lambda(B_{k}) = \sum_{m:{\mathbf w}_{m}\in
  B_{k}}\left(E[Z^{2}_{1,m}]-E[Z^{2}_{2,m}]\right)^{2}.
\label{saito_1&2_class}
\end{equation}
We see that with this $\lambda$, the best basis given by (\ref{best_basis}) is
the basis maximizing the sum of the euclidean distances between the expected
values of all the
basis coordinates for the two classes. 
Now, we observe that the measure of performance (\ref{saito_1&2_class})
of the basis $B_{k}$ does not consider how the data is distributed 
around the expected values. For example, if:
\begin{eqnarray*}
&I_{1,m} =
\left[E[Z^{2}_{1,m}]-\sqrt{Var[Z^{2}_{1,m}]},E[Z^{2}_{1,m}]+\sqrt{Var[Z^{2}_{1,m}]}\right]&
\\
&I_{2,m} =
\left[E[Z^{2}_{2,m}]-\sqrt{Var[Z^{2}_{2,m}]},E[Z^{2}_{2,m}]+\sqrt{Var[Z^{2}_{2,m}]}\right]&
\end{eqnarray*}
where $Var[Z^{2}_{y,m}]$ is empirical variance of $Z^{2}_{y,m}$,
then it may well happen that $I_{1,m}\cap I_{2,m} \neq \emptyset$,
even if ${\mathbf w}_{m}\in \mbox{arg}\max_{B_{k}\in{\mathcal
    D}_{i}\in{\mathcal L}}\lambda(B_{k})$.

Ideally, we want a basis $B$ where the overlap ${\mathcal O}(B)$ given by 
\[ {\mathcal O}(B)=\sum_{m:{\mathbf w}_{m}\in B}\left|I_{1,m}\cap I_{2,m}\right| \] 
is as small as possible. That is a basis which simultanously
is discriminating {\em between} classes and has the opposite property
{\em inside} classes. This motivates the following definition of a new
discrimination measure $\lambda^{\prime}$ by 
\begin{equation} 
  \lambda^{\prime}(B_{k}) = \sum_{m:{\mathbf w}_{m}\in
  B_{k}}\left[\frac{E[Z^{2}_{1,m}]-E[Z^{2}_{2,m}]}
{(Var[Z^{2}_{1,m}]+Var[Z^{2}_{2,m}])^{1/2}}\right]^{2}.
\label{eirik1_1&2_class}
\end{equation}
Note how the performance measure in (\ref{eirik1_1&2_class}) defers
from the measure in (\ref{saito_1&2_class}). 
We see that the numerator in (\ref{eirik1_1&2_class}) measures 
the separability of datapoints {\em between} the classes $1,2$,
and the denominator measures the dispersion of the datapoints
{\em inside} each of the classes $1,2$. Neither of the measures
$\lambda^{\prime},\lambda$ captures differences between classes in sign in the basis
coordinates. To improve on this fact, we define the
measure $\lambda^{\prime\prime}$ by
\begin{eqnarray} 
  \lambda^{\prime\prime}(B_{k})&=&\sum_{m:{\mathbf w}_{m}\in B_{k}}
  \left[ E[(Z_{1,m}-Z_{2,m})^{2}]^{1/2}/\right.\nonumber\\
  &&\left(E[(Z_{1,m}({\mathbf x})-Z_{1,m}({\mathbf
      x}^{\prime}))^{2}]^{1/2}\right.+\nonumber \\
    && \left.\left.E[(Z_{2,m}({\mathbf x}^{\prime\prime})-
    Z_{2,m}({\mathbf x}^{\prime\prime\prime}))^{2}]^{1/2}\right)\right],\nonumber \\
  &&{\mathbf x}\neq{\mathbf x}^{\prime}\in{\mathcal X}^{(1)},\ 
  {\mathbf x}^{\prime\prime}\neq {\mathbf
    x}^{\prime\prime\prime}\in{\mathcal X}^{(2)}.
  \label{eirik2_1&2_class}
\end{eqnarray}
We see that the numerator in (\ref{eirik2_1&2_class}) measures the
separability of {\em signed} datapoints between the classes $1,2$,
and the denominator measures the dispersion of {\em signed} datapoints
{\em inside} these classes.

\subsection{Construction of an Oracle Classifier Using Multiple LDB's}

The construction is due to the following observation: Having chosen a best basis
${\mathbf \Psi}^{t}_{n\times n}$, where ${\mathbf \Psi}_{n\times n} =
\mbox{arg}\max_{B_{k}\in{\mathcal D}_{i}\in{\mathcal L}}\zeta(B_{k})$,
and $\zeta$ is some discrimination measure, there are subsets 
${\mathcal S}_{j}$ of the set ${\mathcal X}$ of input signals on which 
${\mathbf \Psi}^{t}_{n\times n}$ works better than other subsets. That is, the signals
in disjoint sets ${\mathcal S}_{j}$ have significant differences in how they distribute their
energy among the different elements in the basis
${\mathbf \Psi}^{t}_{n\times n}$. More precisely: Let
${\mathcal W}_{K}$ be the {\em feature space} of dimension $K<n$
spanned by the $K$ most important elements in the best basis
${\mathbf \Psi}^{t}_{n\times n}$, sorted in decreasing order of importance, and
$P_{{\mathcal W}_{K}}$ be the orthogonal projection onto  ${\mathcal W}_{K}$.

Now, consider the sets $A^{(k)}$ and $B^{(k)}$ of points in $k$-dimensional
euclidean space given by:
$A^{(k)}=\{P_{{\mathcal W}_{k}}{\mathbf x}_{j}\}_{{\mathbf x}_{j}\in
  {\mathcal X}^{(1)}}\subset [-1,1]^{k}$, 
$B^{(k)}= \{P_{{\mathcal W}_{k}}{\mathbf x}_{j}\}_{{\mathbf x}_{j}\in
  {\mathcal X}^{(2)}}$ $\subset [-1,1]^{k}$. It is clear by the definition
of ${\mathcal W}_{k}$, that the two point-clouds $A^{(k)}$ and $B^{(k)}$ should
be concentrated in more or less disjoint regions in $[-1,1]^{k}$ if the two
classes are separable by our method, that is we should
observe clustering when plotting the points of $A^{(k)}$
and $B^{(k)}$ in $[-1,1]^{k}$ and labeling each point after its class.

We sort out clusters by the following recursive algorithm.
\begin{algthm}
  The Dyadic Cluster Search Algorithm (DCSA). Given appropriately chosen numbers
    $n\geq K\geq 1,\ 1>\delta>0,\ 1> \eta\geq 0,\ 1>\mu\geq\nu>0$.  \\
\item[\bf Step 0:] Choose a performance measure $\lambda$ as in
    (\ref{saito_1&2_class}),(\ref{eirik1_1&2_class}) or
    in (\ref{eirik2_1&2_class}), or some other favourite measure. Set
    $\beta=\lceil\nu\left|{\mathcal X}\right|\rceil,\
    \gamma_{A}=\lceil\eta\left|{\mathcal X}^{(1)}\right|\rceil, \
    \gamma_{B}=\lceil\eta\left|{\mathcal X}^{(2)}\right|\rceil,\
    I=[-1,1]$.  \\
\item[\bf Step 1:] Select the feature spaces ${\mathcal W}_{K}$ by the formula
    (\ref{best_basis}) and truncate to the $K<n$ most important
    basis elements. Compute the sets $A^{(K)},\ B^{(K)}$
    as defined above. Set $\Delta=0.0,\ k=1,\
    C^{(k)}=I^{k},\ C^{(k)}_{next}=I^{k}$, $FoundCluster=0$. \\
\item[\bf Step 2:] Set $A=A^{(k)},\ B=B^{(k)}, \ C=C^{(k)},\ C_{next}=C^{(k)}_{next}$. If
    $\left|A\right|\leq \gamma_{A} \mbox{ and }\left|B\right|\leq
    \gamma_{B}$, terminate the algorithm. Else, compute
    $\alpha=\lceil\mu(\left|A\right|+\left|B\right|)\rceil$, 
    $N_{A}(C)=\left|A\cap C\right|,\ N_{B}(C)=\left|B\cap C\right|$. \\
\item[\bf Step 3:] If $N_{A}+N_{B}\geq \max(\alpha,\beta)$, compute
    the error rate $\epsilon=\min(N_{A},N_{B})/(N_{A}+N_{B})$ and proceed to the next
    step. Else, if $C_{next}\neq C$, set $C=C_{next}$ and jump to {\bf Step 2}. 
    Else, if $C_{next}=C$, and $FoundCluster=1$, jump to {\bf Step 1}.
    Else, if $C_{next}=C$, $FoundCluster=0$, if $k<K$, set $k=k+1$ and jump
    to {\bf Step 2}. Else, if  $C_{next}=C$, $FoundCluster=0$,
    $k=K$, set $k=1$, $\Delta=\Delta+\delta$ and jump to {\bf Step
      2}. \\
\item[\bf Step 4:] If $\epsilon \leq \Delta$, store the location of the cube $C$
    together with the numbers $N_{A},N_{B}$ and identification of the $k$ basis
    elements defining the space ${\mathcal W}_{k}$. 
    Then, for each index $i\in \{1,2,...,K\}$, set
    $A^{(i)}=A^{(i)}-\left(P_{{\mathcal W}_{i}}A^{(K)}\right)\cap C,\
    B^{(i)}=B^{(i)}-\left(P_{{\mathcal W}_{i}}B^{(K)}\right)\cap C$ and
    for each ${\mathbf x}_{j_{i}}\in {\mathcal X}^{(i)}$,
    if $P_{{\mathcal W}_{k}}{\mathbf x}_{j_{i}}\in C,\ {\mathcal X}^{(i)}={\mathcal X}^{(i)}-{\mathbf x}_{j_{i}},\ i=1,2$. Set
    $FoundCluster=1,\ \Delta=0.0,\ k=1$, 
    and jump to {\bf Step 2}. Else, divide $C$ into $2^{k}$ subcubes
    $C_{1},...,C_{2^{k}}$ by splitting 
    each of the sidelengths of $C$ into two sides of equal length, and for
    each index $i=1,2,...,2^{k}$, jump to {\bf Step 2} with $C=C_{i}$,
    $C_{next}=C_{i+1}, 1\leq i<2^{k}, \ C_{next}=C, i=2^{k}$. \\
  \label{DCSA}
\end{algthm}

Less precisely: This algorithm carries out a classification on the signals in the
input signal space ${\mathcal X}={\mathcal X}^{(1)}\cup{\mathcal X}^{(2)}$
by dividing the set ${\mathcal X}$ into disjoint subsets
${\mathcal S}_{j}$ and performing
a classification on each of these subsets represented in a basis
${\mathbf \Psi}_{j}^{t}$. Each
${\mathcal S}_{j}$  consists exclusively of the signals
${\mathbf x}_{i}\in{\mathcal X}$ on which the most discriminating
basis ${\mathbf \Psi}_{j}^{t}$ selected by (\ref{best_basis}) performs
best. Having computed a best basis ${\mathbf \Psi}_{1}^{t}$, the set ${\mathcal S}_{1}$ is
selected first, the signals in ${\mathcal S}_{1}$ are assigned 
class names and then ${\mathcal S}_{1}$ is deleted from the set
${\mathcal X}$. 
Then a new best basis ${\mathbf \Psi}_{2}^{t}$ for the new ${\mathcal X}$ is computed by the 
formula (\ref{best_basis}), the set ${\mathcal S}_{2}$ is selected, and
so on. The algorithm terminates when 
the set ${\mathcal X}$ has become sparse. Thus, we see that by adapting the parameters
we can prevent the algorithm from trying to classify the 
part of the training dataset which it finds most difficult to
classify, and so we gain a smaller overall
training-error-rate. But this adjusting of parameters has to be done
carefully, so that the algorithm does not fail to catch important
features of the signals. The algorithm selects the
subsets ${\mathcal S}_{j}$ using as few features as possible, starting
with only the most important feature element (= the most discriminating
basis element in the best basis). Then, given some upper limit on the rate of
error allowed in the clusters, if no clean clustering is observed in
the feature space of this single feature element, the algorithm adds information by
taking into consideration also the second best feature element and
looks for clustering in the feature space spanned by the two best feature
elements and so on. If no clean clustering is observed using all $K$ best feature
elements, the upper error limit is increased and the feature space of
the one most important feature element is again searched for clusters,
and so on. Using as few features as possible reduces the risk of 
overtraining of the algorithm, that is the algorithm selecting features 
that are too adapted to the specific set ${\mathcal X}$ of training
data. On the other hand, we see that this algorithm is flexible
in its selection of relevant features in that it constructs a
sequence of feature extractors $\{{\mathbf \Psi}^{t}_{j}\}$ where each
${\mathbf \Psi}^{t}_{j}$ is
specially adapted to some part ${\mathcal S}_{j}$ of the dataset
${\mathcal X}$. The output of the algorithm is a sequence  ${\mathcal
  C}=\{C_{i}\}_{i=1}^{L}$ of dyadic hypercubes $C_{i}\subset[-1,1]^{K}$ of possibly different
dimensions $k_{i}, 1\leq k_{i}\leq K,\ 1\leq i \leq L$, where to each
cube $C_{i}$ corresponds a specific feature space ${\mathcal W}_{k_{i}}$ as defined
above, and a class name $y_{C_{i}}$ which equals the name of the majority
class of the set
$\left(\{P_{{\mathcal W}_{k_{i}}}{\mathbf x}_{j}\}_{{\mathbf
      x}_{j}\in{\mathcal X}}\right)\cap
C_{i}=\{P_{{\mathcal W}_{k_{i}}}{\mathbf x}_{j}\}_{{\mathbf
    x}_{j}\in{\mathcal S}_{i}}$
of datapoints in ${\mathcal W}_{k_{i}}$ that $C_{i}$ contains. 
We will call ${\mathcal C}$ a simple two-class {\em oracle classifier},
or simply oracle, for the two-class problem
$d:{\mathcal X}^{(1)}\cup{\mathcal X}^{(2)}\rightarrow
{\mathcal Y}=\{1,2\}$. 

\subsection{On Using and Choosing Oracle Classifiers}

Given a two-class problem
$d:{\mathcal X}^{(1)}\cup{\mathcal X}^{(2)}\rightarrow
{\mathcal Y}=\{1,2\}$, we compute ${\mathcal C}=\{C_{j}\}_{j=1}^{L}$
by the DCSA. Then, given a sample ${\mathbf x}\in{\mathcal T}$, where ${\mathcal T}$ is a
test dataset, we assign ${\mathbf x}$ to a class by the following procedure: 
We check if: $P_{{\mathcal W}_{k_{j}}}{\mathbf x}\in C_{j}$, starting with index
$j=1$ and continuing until we get a positive
answer for some index $j^{\prime}\leq L$. We then assign a weighted class
$y_{C_{j^{\prime}}}$-vote to ${\mathbf x}$ by
computing the product of $1-\epsilon_{j^{\prime}}$, where
$\epsilon_{j^{\prime}}$ is the error rate of $C_{j^{\prime}}$, and its
statistical frequency
$(N_{A}(C_{j^{\prime}})+N_{B}(C_{j^{\prime}}))/\left|{\mathcal X}\right|$.  
If $P_{{\mathcal W}_{k_{j}}}{\mathbf x}\not\in C_{j},\ \forall
C_{j}\in{\mathcal C}$, we consider the class of ${\mathbf x}$ undetermined.

Different choices of discrimination measure or different
settings of the parameters in the DCSA result in different
classifiers. For a two-class problem, we
can construct several classifiers by using different performance measures/parameters,
and let the weighted majority vote of the classifiers decide whether a sample
${\mathbf x}\in{\mathcal T}$ is of class 1 or class 2. For a $n$-class problem, $n>2$,
we will apply the method of splitting the $n$-class problem into $n$ two-class problems:
$d:{\mathcal X}\rightarrow \{i,0\},\ 1\leq i \leq n$,
as proposed in Saito, Coifman (1996), by splitting the training data set into
two sets of class $i$ and not $i$. One then constructs oracles for each
two-class problem. To classify an unknown sample
${\mathbf x}\in{\mathcal T}$, we compute weighted class
votes as explained above for the set
of oracles and assign ${\mathbf x}$ to the majority vote class.

\section{Experimental Results}
 
In some of the calls to the DCSA in the experiments described below we allowed
the algorithm to select a best basis only once, we call this method a 
LDB-method (Local Discriminant Basis-method). 
In the cases were we allowed the algorithm to select multiple
different best bases in sequence, we call the method a MLDB-method
(Multiple Local Discriminant Basis-method). In the cases
where we organized the classifiers resulting from different calls (calls
with different discrimination measures) to the DSCA into a classifier
by taking the majority vote over these classifiers, we call the method
a {\em superposition} LDB or MLDB-method, denoted SLDB or
SMLDB-method, respectively. In all the three examples below we
generated 10 independent realizations of both the training dataset and the
test dataset. The results shown in Table \ref{table1}, Table
\ref{table2}, Table \ref{table3} are the mean over the 10 simulations corresponding to the
10 independent realizations of the datasets. 

\subsection{Example 1} 

We consider a two class waveform classification problem as presented 
in Fossgaard (1997).
We generated sets of 100 training signals and 1000 test signals of length 1024
for each class by the formula 

\begin{eqnarray}
 && Q_{n}(R,\theta,t) = \nonumber \\
  && C(R,t)\sum_{j=1}^{n}A_{n}(j)
  e^{ik(\frac{r_{j}^{2}}{2R}-r_{j}\cos(\theta-\theta_{j}))},
  \label{approximated_signal}
\end{eqnarray}
\noindent where we have:
\begin{eqnarray*}
&& C(R,t)=\frac{e^{-ik(ct-R)}}{R}\mbox{ is considered constant $=1$}\\
&& \mbox{for simplicity. } \\
&& R = 10^{4}.  \\
&& k = 100. \\
&& A_{n}(j) = \frac{1}{n}. \\
&& r_{j} \mbox{ is random variable uniformly distributed on } [1,10]. \\
&& \theta_{j} \mbox{ is random variable uniformly distributed on}\\
&& [2\pi\frac{j}{n},2\pi\frac{j}{n} +\frac{\pi}{4}]. 
\end{eqnarray*}
For each $n$-tiple of realizations $\{r_{j},\theta_{j}\}_{j=1}^{n}$ of
the pair of random variables $r_{j},\theta_{j}$, we generate a discrete signal
$S_{n}(\theta)$ by uniformly sampling the real part of $Q_{n}(R,\theta,t)$ 1024
times in the variable $\theta$ with sampling density $2\pi/16\cdot k = 2\pi/1600$.
We generated data sets by extracting realizations of $S_{n}/\|S_{n}\|_{2}$
{\em smoothly} from a fixed sampling interval. In this problem we used
$n=3,\ n=4$ in (\ref{approximated_signal}) to define two classes of
signals and the coiflet with filterlength 18 as dictionary.  All calls to the DCSA
in this experiment were made with $K=5,\ \delta = 0.01,\ \eta = 
0.05, \ \mu = 0.10, \ \nu = 0.05$.  The results are shown in Table \ref{table1}.

\begin{table}[h]
\begin{tabular}{lcccccccc} \hline\hline
   Method & 
  \multicolumn{4}{c}{ Classification rate (\%)} &
  \multicolumn{4}{c}{ Error rate (\%)} \\
  & \multicolumn{2}{c}{ Training data} &
  \multicolumn{2}{c}{ Test data} 
  & \multicolumn{2}{c}{ Training data} &
  \multicolumn{2}{c}{ Test data} \\ 
  &  Total &  $\sigma$ & 
   Total &  $\sigma$ & 
   Total &  $\sigma$ 
  &  Total &  $\sigma$ \\ \hline
 LDB1  &99.7 &0.7 &99.5 &0.9 &19.9 &3.0 &29.8 &2.8  \\ \hline
 MLDB1 &98.5 &1.4 &98.4 &1.4 &8.6 &1.1 &23.5 &2.7   \\ \hline
 LDB2  &97.6 &2.0 &96.3 &3.6 &16.0 &4.5 &23.5 &4.4   \\ \hline
 MLDB2 &95.2 &1.9 &93.8 &2.3 &12.9 &3.1 &23.7 &3.4  \\ \hline
 LDB3  &98.9 &1.7 &98.8 &2.1 &16.6 &4.6 &24.6 &4.9  \\ \hline
 MLDB3 &98.4 &1.1 &99.5 &0.6 &13.7 &3.1 &24.4 &1.9   \\ \hline
 SLDB  &100 &0.0 &100 &0.0 &14.7 &3.9 &22.2 &1.5   \\ \hline 
 SMLDB &100 &0.0 &100 &0.0 &9.1 &3.3 &{\bf 20.4} &2.0  \\ \hline\hline
\end{tabular}
\caption{ The average classification rates and the
  corresponding error rates over 10 simulations from Example 1. 
  LDB1 is the LDB selected by the measure
  $\lambda^{\prime}$. MLDB1 is the MLDB selected by  
  the measure $\lambda^{\prime}$. LDB2 is the LDB selected by the 
  measure $\lambda^{\prime\prime}$. MLDB2 is the MLDB selected by the  
  measure $\lambda^{\prime\prime}$. LDB3  is the LDB selected by the 
  measure  $\lambda$. MLDB3 is the MLDB selected by the measure $\lambda$.
  SLDB is the superposition of methods LDB1, LDB2, LDB3. SMLDB is the superposition
  of the methods MLDB1, MLDB2, MLDB3. $\sigma$ is the square root of
  the sample variance. }
\label{table1}
\end{table}

\subsection{Example 2} 

This example is identical to Example 1 except that we used $n=4,\ n=5$
in (\ref{approximated_signal}) to define the two signal classes. We
used the coiflet with filterlength 18 as dictionary. All calls to the
DCSA in this experiment were made with $K=5,\ \delta = 0.01,\ \eta = 
0.05, \ \mu = 0.10, \ \nu = 0.05$.  The results are shown in Table
\ref{table2}.

\begin{table}[h]
 \begin{tabular}{lcccccccc} \hline\hline
   Method & 
  \multicolumn{4}{c}{ Classification rate (\%)} &
  \multicolumn{4}{c}{ Error rate (\%)} \\ 
  & \multicolumn{2}{c}{ Training data} &
  \multicolumn{2}{c}{ Test data} 
  & \multicolumn{2}{c}{ Training data} &
  \multicolumn{2}{c}{ Test data} \\ 
  &  Total &  $\sigma$ & 
   Total &  $\sigma$ & 
   Total &  $\sigma$ 
  &  Total &  $\sigma$ \\ \hline
 LDB1  &99.3 &1.4 &98.7 &2.3 &11.9 &3.5 &19.9 &5.4  \\ \hline
 MLDB1 &98.0 &1.4 &97.3 &2.3 &6.4 &2.9 &20.5 &4.3   \\ \hline
 LDB2  &96.9 &2.6 &96.9 &2.5 &10.5 &2.7 &{\bf 17.5} &3.2   \\ \hline
 MLDB2 &96.0 &2.2 &94.1 &3.6 &9.5 &1.9 &19.0 &2.6  \\ \hline
 LDB3  &99.1 &1.7 &99.6 &0.9 &24.2 &5.6 &32.6 &7.9  \\ \hline
 MLDB3 &98.5 &1.4 &99.5 &0.5 &21.6 &3.6 &35.8 &4.2   \\ \hline
 SLDB  &100 &0.0 &100 &0.0 &15.0 &4.9 &21.8 &4.5   \\ \hline 
 SMLDB &100 &0.0 &100 &0.0 &8.4 &3.7 &20.1 &3.2  \\ \hline\hline
\end{tabular}
\caption{    The average classification rates and the
  corresponding error rates over 10 simulations from Example 2. 
  LDB1 is the LDB selected by the measure
  $\lambda^{\prime}$. MLDB1 is the MLDB selected by  
  the measure $\lambda^{\prime}$. LDB2 is the LDB selected by the 
  measure $\lambda^{\prime\prime}$. MLDB2 is the MLDB selected by the  
  measure $\lambda^{\prime\prime}$. LDB3  is the LDB selected by the 
  measure  $\lambda$. MLDB3 is the MLDB selected by the measure $\lambda$.
  SLDB is the superposition of methods LDB1, LDB2, LDB3. SMLDB is the superposition
  of the methods MLDB1, MLDB2, MLDB3. $\sigma$ is the square root of
  the sample variance. }
\label{table2}
\end{table}

\subsection{Example 3}

We consider a three class waveform classification problem as
presented in Saito (1994). 
We generated sets of 100 training signals and 1000 test signals of
length 32 for each class by first extracting signal samples by the formulas 
\begin{eqnarray*} 
  && f_{1}(i) = uh_{1}(i)+(1-u)h_{2}(i)+\epsilon(i) \mbox{ for Class
    1}  \\
  && f_{2}(i) = uh_{1}(i)+(1-u)h_{3}(i)+\epsilon(i) \mbox{ for Class
    2}  \\
  && f_{3}(i) = uh_{2}(i)+(1-u)h_{3}(i)+\epsilon(i) \mbox{ for Class
    3}, 
\end{eqnarray*}
where $i=1,...,32,\ h_{1}(i)=\max(6-|i-7|,0),\
h_{2}(i)=h_{1}(i-8), \ h_{3}(i)=h_{1}(i-4),$ $u$ is a uniform random
variable on the interval $(0,1)$, and $\epsilon(i)$ are the standard
normal variates. We then normalized the signals in the energy norm by
setting $f_{1}(i)=f_{1}(i)/\|f_{1}\|_{2},\ f_{2}(i)=f_{2}(i)/\|f_{2}\|_{2},\
f_{3}(i)=f_{3}(i)/\|f_{3}\|_{2},\ i=1,...,32$.  
We used the coiflet with filterlength 6 as a dictionary for this
problem. All calls to the DCSA in this experiment were made with $K=5,\ \delta = 0.01,\ \eta =
0.05, \ \mu = 0.20, \ \nu = 0.05$. The results are shown in Table \ref{table3}. 
%

\begin{table}[h]
 \begin{tabular}{lcccccccc} \hline\hline
   Method & 
  \multicolumn{4}{c}{ Classification rate (\%)} &
  \multicolumn{4}{c}{ Error rate (\%)} \\
  & \multicolumn{2}{c}{ Training data} &
  \multicolumn{2}{c}{ Test data} 
  & \multicolumn{2}{c}{ Training data} &
  \multicolumn{2}{c}{ Test data} \\ 
  &  Total &  $\sigma$ & 
   Total &  $\sigma$ & 
   Total &   $\sigma$ &
   Total &   $\sigma$ \\ \hline
 LDB1  &100 &0.0 &100 &0.0 &23.7 &1.9 &28.2 &0.8  \\ \hline
 MLDB1 &100 &0.0 &100 &0.0 &22.9 &2.1 &28.5 &2.7  \\ \hline
 LDB2  &100 &0.0 &100 &0.0 &23.6 &2.4 &27.9 &1.9 \\ \hline
 MLDB2 &100 &0.0 &100 &0.0 &20.7 &2.7 &27.7 &1.5  \\ \hline
 LDB3  &100 &0.0 &100 &0.0 &25.0 &2.6 &29.1 &1.8  \\ \hline
 MLDB3 &100 &0.0 &100 &0.0 &23.0 &2.6 &26.2 &2.6   \\ \hline
 SLDB  &100 &0.0 &100 &0.0 &18.7 &1.9 &22.7 &0.9   \\ \hline 
 SMLDB &100 &0.0 &100 &0.0 &15.7 &1.9 &{\bf 20.5} &1.0 \\ \hline\hline
\end{tabular}
 \caption{   The average classification rates and the
    corresponding error rates over 10
    simulations from Example 3. LDB1 is the LDB selected by the measure
    $\lambda^{\prime}$. MLDB1 is the MLDB selected by the measure 
    $\lambda^{\prime}$. LDB2 is the LDB selected by the 
    measure $\lambda^{\prime\prime}$. MLDB2 is the MLDB selected by the
    measure $\lambda^{\prime\prime}$. LDB3 is the LDB selected by the 
    measure  $\lambda$. MLDB3 is the MLDB selected by the measure $\lambda$.
    SLDB is the result from a superposition of the methods LDB1, LDB2,
    LDB3. SMLDB is the result from a superposition of the methods
    MLDB1, MLDB2, MLDB3.
    $\sigma$ is the square root of the sample variance.}
\label{table3}
\end{table}

\section{Comments}

\subsection{Comments to Example 1}

In this example we achieved the best result by the superposition method
using multiple LDB's, denoted SMLDB. We see that the generalized methods MLDB1,
MLDB2, MLDB3 are almost indistinguishable in this example, we conclude
that our new measures $\lambda^{\prime},\ \lambda^{\prime\prime}$ 
hardly yield a significantly better classification than the original measure
$\lambda$, the positive effect is in any case small. 
Furtermore, for the measure $\lambda^{\prime}$ we do get better
results by the generalized method, whereas for the measures
$\lambda^{\prime\prime}$ and $\lambda$ the positive effect of
generalizing is more doubtful. But all in all, it seems we are a
little better off with either measure $\lambda^{\prime},\
\lambda^{\prime\prime}$ than the original $\lambda$.

\subsection{Comments to Example 2}

In this example we achieved the best result with the method LDB2.
We see that both discrimination measures $\lambda^{\prime},\
\lambda^{\prime\prime}$ clearly outperform the original
measure $\lambda$ in this problem. As in the previous example, the
measures $\lambda^{\prime}$ and $\lambda^{\prime\prime}$ yield about the same
results with MLDB-methods. When not taking superpositions of several
classifiers, the generalised MLDB-method does not yield any
improvements in results on test data, rather it seems that this method
adapts too much to training data in this example. Furthermore, due
to the poor performance of the measure $\lambda$ in this example, we
get worse results with superposition methods in this example than
when using the best single classifier. But we could expect to further
lower the best error rate on test data by combining classifiers from the
measures $\lambda^{\prime},\ \lambda^{\prime\prime}$ only.

\subsection{Comments to Example 3}

In this example we achieved the best result by the method SMLDB, and
we see that superposition methods are clearly favourable in this
case. However, it seems to make little difference which measure we are using
when not taking superpositions of several classifiers. We remark that
both the measures $\lambda^{\prime},\ \lambda^{\prime\prime}$ select
the standard basis as the most discriminating basis in the first steps
in the DCSA, whereas $\lambda$ does not choose this basis in any step.

\subsection{Conclusion}

We have shown that estimating expectations and variances directly from
the expansion coefficients of the datasets in the binary-tree
structured dictionary of bases may lead to better results than when
using the energy-density dictionaries of bases. Also, we have shown 
that comparing/combining different discrimination measures in
classification problems may lead to significant improvements in the 
success of the classification methods.



%



%











\appendix

\section{Applied software and hardware}

All algorithms and transforms used in the numerical experiments, except
some of the random number generators described below, were implemented in the computer
language C$++$ and compiled with the GNU project C$++$ compiler on a HP
K260 machine with a PA 8000 processor.

\subsection{Random number generators}

In the examples we used the Fortran NAG-routines
G05DAF, G05FAF for generating random numbers with uniform distribution, and
G05FDF for generating random numbers with standard normal
distribution.



































\end{document}